\documentclass[a4paper,reqno]{amsart}

\usepackage{amsthm}
\usepackage{amsmath}
\usepackage{amssymb}
\usepackage{mathrsfs}
\usepackage{latexsym}
\usepackage{exscale}
\usepackage{geometry} % margenes
\usepackage{graphicx} % graficos
\usepackage[utf8]{inputenc}
\usepackage{verbatim}
\usepackage{enumerate} % enumerar con letras
\usepackage{fancyhdr}
\usepackage{multicol}
\usepackage{bbm}% escribir el 1  con mathbb

\usepackage{color}
\definecolor{red}{rgb}{1,0.1,0.1}%valores de las componentes roja, verde y azul
\definecolor{blue}{rgb}{0.1,0.1,1}
\definecolor{vb}{RGB}{160,32,240}
\theoremstyle{plain}

\newtheorem*{teo*}{Theorem}
\newtheorem*{prop*}{Proposition}
\newtheorem*{lemma*}{Lemma}

\numberwithin{equation}{section}
\newtheorem{teo}{Theorem}[section]
\newtheorem{lemma}[teo]{Lemma}
\newtheorem{corol}[teo]{Corollary}
\newtheorem{prop}[teo]{Proposition}

\theoremstyle{remark}
\newtheorem{remark}[teo]{Remark}

\theoremstyle{definition}

\newtheorem*{mydef*}{Definition}

\calclayout

\allowdisplaybreaks

\title[Fractional Fefferman--Stein inequalities on trees]
{Fractional Fefferman--Stein inequalities on infinite rooted $k$-ary trees}

\author{Gonzalo Ibañez-Firnkorn}
 \address{G.~Iba\~{n}ez~Firnkorn. 
Instituto de Matemática (INMABB), Departamento de Matemática, Universidad Nacional del Sur (UNS)-CONICET
\\Bahía Blanca, Argentina}
\email{gonzalo.ibanez@uns.edu.ar}
\author{Jorgelina Recchi}
\address{Jorgelina Recchi. 
Instituto de Matemática (INMABB), Departamento de Matemática, Universidad Nacional del Sur (UNS)-CONICET
\\Bahía Blanca, Argentina}
\email{drecchi@uns.edu.ar}


\begin{document}

\begin{abstract}
We study weighted endpoint estimates for the centered fractional maximal operator on the
infinite rooted $k$-ary tree, $T$. Our first result is a Fefferman--Stein type estimate
\[
w\bigl(\{M_\alpha f>\lambda\}\bigr)
\lesssim_{\alpha,s}
\frac1\lambda\sum_{x\in T}|f(x)|M_{\alpha,s}w(x),
\qquad s>1,
\]
with a constant independent of the branching parameter $k$. Although the proof is
inspired by the non-fractional argument of Ombrosi, Rivera-R\'ios and Safe in \cite{ORRS}, the
fractional normalization introduces a genuinely new scale interaction. At radius $r$,
one must renormalize the function by $k^{r\alpha}$ while keeping the dyadic cutoff at
$k^r$; this leads to a new summation mechanism involving simultaneously
$k^{r\alpha}$ and $k^{r(1-\alpha)}$.

We then study sharpness and the geometry of admissible weights. We construct
counterexamples showing that the estimate cannot in general be reduced to $s=1$
for any $0<\alpha<1$. In the subcritical range $0<\alpha<1/2$, we characterize
power-type radial weights satisfying the natural fractional $A_1$ condition and
prove that the failure persists after any fixed finite number of iterations of
$M_\alpha$. We also prove that
$M_\alpha w\lesssim w$ is necessary for the corresponding weighted weak type estimate
but is not sufficient in the subcritical range. Finally, we identify a phase transition:
for $1/2\le\alpha<1$ no nonzero global weight can satisfy $M_\alpha w\lesssim w$.
\end{abstract}

\maketitle

\section{Introduction and main results}

 Let $k\geq2$ be an integer. We consider $T$ the infinite rooted $k$-ary tree, that is, a graph with an initial vertex $x_0$, called root, and each vertex has $k$ children.

Throughout the paper, we denote $A\lesssim_{\theta} B$ to mean that there exists a positive constant $c_{\theta}$ depending only on $\theta$ such that $A\leq c_{\theta} B$. We say that $A\simeq_{\theta} B$ when $A\lesssim_{\theta} B$ and $B\lesssim_{\theta} A$. 

It is possible to equip $T$, with the graph distance $d$ and the counting measure $\mu$. The distance $d(x,y)$ is the number of edges of the unique path between $x$ and $y$ and the measure $\mu$ counts the number of vertices in the set. If $A\subset T$, we denote $|A|=\mu(A)$ and 
\[f(A)=\int_A f(x) \,dx = \sum_{x\in A} f(x).\]

All radii are nonnegative integers. By a weight we mean a function
$w:T\to[0,\infty)$ which is finite at every vertex; it may vanish on a subset of
$T$. We write $w(E)=\sum_{x\in E}w(x)$.

Let $S(x,r)$ denote the sphere of radius $r$, $S(x,r)=\{y\in T : d(x,y)=r\}$, observe that $S(x,r)$ has exponential growth,
\[
 |S(x,r)|\simeq k^r,
\]
 and $\mu$ is not a doubling measure. This geometry makes $T$ a natural
model in which classical covering arguments are unavailable but sharp maximal estimates
can still be obtained by exploiting combinatorial expansion.

The study of classical operators in the context of $k$-tree has recently received considerable attention, particularly through the works of Naor and Tao \cite{NT10} and  Ombrosi, Rivera-R\'ios and Safe \cite{ORRS}. We also refer to the articles \cite{ST16,ST19} in which the authors studied connections between the geometrical properties of finite and infinite graphs and the boundedness of the Hardy-Littlewood maximal operator. For weighted estimates on trees, we refer to \cite{GR23,LS24,ORR,YYZ}.

In \cite{NT10}, Naor and Tao proved the unweighted weak type $(1,1)$ for the Hardy-Littlewood maximal operator in the $k$-tree $T$,
\[
 \bigl|\{x\in T: Mf(x)>\lambda\}\bigr|
 \leq \frac{c}{\lambda} \sum_{x\in T}|f(x)|
\]
where $M$ is defined by 
\[
 M f(x)
 =\sup_{r\ge0}\frac1{|B(x,r)|}
   \sum_{y\in B(x,r)}|f(y)|.
\]

Recently, in \cite{ORRS} Ombrosi, Rivera-R\'ios and Safe  proved a weighted
Fefferman--Stein inequality for the Hardy--Littlewood maximal operator on $T$ of the
form
\[
 w\bigl(\{Mf>\lambda\}\bigr)
 \lesssim_s \frac1\lambda\sum_{x\in T}|f(x)|M_s w(x),
 \qquad s>1,
\]
where $M_s(w)=M(w^s)^{1/s}$. The authors showed that not only it is not possible to choose $s=1$ but also it is not even possible to choose any number of iterations of the maximal function for the inequality to hold. Their proof
combines a weighted incidence estimate on spheres with a dyadic decomposition whose
crucial gain is a positive power of $2^n/k^r$.

 The Fefferman-Stein estimate was originally proved in \cite{FS71}, later these estimates and related variants have been studied in a variety of contexts \cite{ABN11,LP14,OP02,P19,T04}.
In the fractional setting, the Fefferman-Stein for the fractional maximal operator is a consequence of classical covering lemmas \cite{FS71}. In \cite{CPSS05}, the authors proved that the Fefferman-Stein estimate for the fractional integral holds with a fractional $L\log L$-maximal operator and cannot be replaced with the fractional maximal operator.

The purpose of the present paper is to develop the corresponding theory for the
fractional maximal operator with $0<\alpha<1$, 
\[
 M_\alpha f(x) =\sup_{r\ge0}\frac1{|B(x,r)|^{1-\alpha}}
   \sum_{y\in B(x,r)}|f(y)|,
\]
and its spherical counterpart 
\[
 M_\alpha^o f(x)
 =\sup_{r\ge0}\frac1{|S(x,r)|^{1-\alpha}}
  \sum_{y\in S(x,r)}|f(y)|.
\]
 Note that, in contrast to the standard Euclidean setting, here it makes sense to consider this kind of maximal function because $S(x,r)$ are not sets of measure $0$. For $k\geq2$, we have $M_\alpha \simeq_\alpha M_\alpha^o$ as we will show in   Proposition \ref{Malpho}. 

The fractional problem is not obtained by a
formal substitution in the proof of \cite{ORRS}. The factor $|S(x,r)|^\alpha\simeq
k^{r\alpha}$ changes in a substantial form the distributional decomposition and the final summation.

\medskip

We consider the fractional $L^s$ maximal operator
\[
 M_{\alpha,s}w(x)
 =\sup_{r\ge0}|B(x,r)|^\alpha
 \left(\frac1{|B(x,r)|}\sum_{y\in B(x,r)}w(y)^s\right)^{1/s}.
\]

 The first result gives a
  weighted endpoint estimate  that survives fractional normalization in a space with exponential volume growth.  The following theorem provides a Fefferman-Stein estimate with $M_{\alpha,s}w$ in the right-side.
\begin{teo}\label{teo fS con s}
Let $k\ge2$, $0<\alpha<1$ and $s>1$. Then, for every weight $w\ge0$ on $T$,
\begin{equation} \label{eq:FS-main}
 w\bigl(\{x\in T:M_\alpha f(x)>\lambda\}\bigr)
 \le C_{\alpha,s}\frac1\lambda
 \sum_{x\in T}|f(x)|M_{\alpha,s}w(x),
\end{equation}
where $C_{\alpha,s}$ is independent of $k$. The constant obtained in the proof
deteriorates as $s\downarrow1$.
\end{teo}

Theorem \ref{teo fS con s} is proved in Section~2. The weighted sphere-incidence lemma, Lemma \ref{1otimesw}, is a fractional
variant of the combinatorial estimate in \cite{ORRS}, but its power of $k^r$ changes from
the non-fractional one. 

A central issue is whether the assumption $s>1$ in Theorem \ref{teo fS con s} is genuinely necessary. The following counterexample therefore plays a decisive role: in the subcritical range $0<\alpha<1/2$, it shows that the estimate fails not only for $s=1$, but even when $M_\alpha$ is replaced by any fixed finite iterate. This establishes the necessity of the modification introduced in Theorem \ref{teo fS con s} and identifies a fundamental obstruction in the fractional setting.

\begin{teo} \label{Teo 1.2}
For any $n\in \mathbb{N}$ and $0<\alpha<1/2$, there exists a weight $w$ such that $M_\alpha w\lesssim w$ and 
a sequence of finitely supported functions $\{f_j\}$ such that
\[
 \frac{w(\{M_\alpha f_j>1\})}
 {\displaystyle\sum_{x\in T}|f_j(x)|M_\alpha^{(n)}w(x)}
 \xrightarrow[j \to \infty]{} \infty.
\]

Consequently, no estimate of the form
\[
w(\{M_\alpha f>\lambda\})
\leq \frac{C}{\lambda}\sum_{x\in T}|f(x)|M_\alpha^{(n)}w(x)
\]
can hold uniformly for all finitely supported $f$ and all $\lambda>0$.

\end{teo}

Although Theorem \ref{Teo 1.2} is the fractional counterpart of the result of \cite{ORRS}, the construction of the counterexample requires a different approach. In their example, the weight is chosen to compensate the average of the ball. 
In the fractional setting, however, this normalization becomes $|B(x,r)|^{1-\alpha}$, and the natural adaptations of their weight no longer yield the balance needed to produce a coun\-ter\-exam\-ple. Thus, the non-fractional construction does not extend directly to our setting. Overcoming this difficulty requires a more delicate example, based on a suitably chosen subtree whose geometry compensates for the fractional normalization.

We next obtain a necessary condition for the weighted weak-type $(1,1)$ estimate.
Testing this estimate on point masses recovers the condition $M_\alpha w\lesssim w$.
The preceding counterexample shows that this condition is not sufficient in the
subcritical range.

\begin{teo}\label{Teo 1.3}
Let $0<\alpha<1$ and let $w$ be a weight on $T$. Assume that there exists
$C_w>0$ such that
\begin{equation}\label{eq:weak-w-assumption}
w\bigl(\{x\in T:M_\alpha f(x)>\lambda\}\bigr)
\leq \frac{C_w}{\lambda}\sum_{x\in T}|f(x)|w(x)
\end{equation}
for every finitely supported $f$ and every $\lambda>0$. Then
\[
M_\alpha w(x)\lesssim_{\alpha,k,C_w}w(x),
\qquad x\in T.
\]
\end{teo}

Consequently, for $0<\alpha<1/2$ this condition is necessary but
not sufficient. If $1/2\le\alpha<1$, no nonzero weight can satisfy
$M_\alpha w\lesssim w$, and hence no nontrivial weighted weak type estimate of the
above form is possible; see Proposition~\ref{prop:no-global-weight}.

At the critical value $\alpha=1/2$, the exponential rates used in the
supercritical argument coincide. The radial computation below nevertheless
provides an additional factor depending on the generation, which rules out
global weights also at the endpoint, without any radiality assumption.

We present the power-type radial weights of the form 
\[
 w_\beta(x)=k^{\beta|x|},
\]
where $|x|$ is the generation of $x$ and $\beta \in \mathbb{R}$.

These weights allow us to compare the necessary condition with the sufficient
condition involving $M_{\alpha,s}$. Contributions from descendants and from
paths passing through ancestors impose two restrictions on $\beta$. Their balance
determines the admissible interval and explains the threshold $\alpha=1/2$.

\begin{teo} \label{teo B}
If $0<\alpha<1/2$, then
\[
 M_\alpha w_\beta\simeq w_\beta
 \quad\Longleftrightarrow\quad
 \beta\in[\alpha-1,-\alpha].
\]
If $1/2\le\alpha<1$, no power-type radial weight satisfies this estimate globally.
Moreover, for $s>1$,
\[
 M_{\alpha,s}w_\beta\simeq w_\beta
\]
holds precisely in the range
\[
 s<\frac1{2\alpha},
 \qquad
 \beta\in\left[\alpha-\frac1s,-\alpha\right].
\]
\end{teo}

The paper is organized as follows. Section~2 contains the comparison between ball and
spherical fractional maximal operators, the weighted incidence estimate, the new
fractional distributional lemma, and the proof of Theorem \ref{teo fS con s}. Section~3 is devoted to
radial weights, sharpness, counterexamples, necessary conditions, and the phase
transition at $\alpha=1/2$.

\section{Proof of Theorem \ref{teo fS con s}}

To prove Theorem \ref{teo fS con s} we need several weighted estimates.
We first compare the ball and spherical fractional maximal operators.
Spheres are convenient for the argument because pairs of vertices at a fixed
distance can be classified by their common ancestor. We keep track of the
constants in this comparison so that the final estimate is uniform in $k$.

We write
\[
 B(x,r)=\{y\in T:d(x,y)\le r\},
 \qquad
 S(x,r)=\{y\in T:d(x,y)=r\}.
\]
The geometry of the rooted tree gives, uniformly in $x$, $r$ and $k\ge2$,
\begin{equation}\label{eq:sphere-size}
 k^r\le |S(x,r)|\le2k^r,
 \qquad
 |B(x,r)|\le2|S(x,r)|.
\end{equation}

For $0<\alpha<1$ and $1\le s<\infty$, we set
\[
 M_{\alpha,s}f(x)
 =\sup_{r\ge0}|B(x,r)|^\alpha
 \left(\frac1{|B(x,r)|}\sum_{y\in B(x,r)}|f(y)|^s\right)^{1/s}
\]
and
\[
 M_{\alpha,s}^o f(x)
 =\sup_{r\ge0}|S(x,r)|^\alpha
 \left(\frac1{|S(x,r)|}\sum_{y\in S(x,r)}|f(y)|^s\right)^{1/s}.
\]
Thus $M_{\alpha,1}=M_\alpha$ and $M_{\alpha,1}^o=M_\alpha^o$.

\begin{prop}\label{Malpho}
Let $0<\alpha<1$. Then
\[
 M_\alpha f(x)\le C_\alpha M_\alpha^o f(x),
 \qquad
 M_\alpha^o f(x)\le2^{1-\alpha}M_\alpha f(x),
\]
where $C_\alpha$ is independent of $k$. Moreover, for every $1\le s<\infty$,
\begin{equation}\label{eq:Mas-sphere-ball}
 M_{\alpha,s}^o f(x)
 \le 2^{(1/s-\alpha)_+}M_{\alpha,s}f(x).
\end{equation}
If, in addition, $\alpha s<1$, then
\begin{equation}\label{eq:Mas-ball-sphere}
 M_{\alpha,s}f(x)\le C_{\alpha,s}M_{\alpha,s}^o f(x),
\end{equation}
where $C_{\alpha,s}$ is independent of $k$.
\end{prop}

\begin{proof}
 By \eqref{eq:sphere-size}, for every $r\ge0$,
\begin{align*}
\frac1{|S(x,r)|^{1-\alpha}}\sum_{S(x,r)}|f|
&\le
\left(\frac{|B(x,r)|}{|S(x,r)|}\right)^{1-\alpha}
\frac1{|B(x,r)|^{1-\alpha}}\sum_{B(x,r)}|f|
\le2^{1-\alpha}M_\alpha f(x).
\end{align*}
Hence $M_\alpha^o f(x)\le2^{1-\alpha}M_\alpha f(x)$.

Conversely,
\begin{align*}
\frac1{|B(x,r)|^{1-\alpha}}\sum_{B(x,r)}|f|
&=\frac1{|B(x,r)|^{1-\alpha}}
  \sum_{\ell=0}^r\sum_{S(x,\ell)}|f|\\
&\le M_\alpha^o f(x)
 \frac{\sum_{\ell=0}^r|S(x,\ell)|^{1-\alpha}}{|B(x,r)|^{1-\alpha}}\\
&\le M_\alpha^o f(x)
 2^{1-\alpha}\sum_{j=0}^r k^{-j(1-\alpha)}\\
&\le \frac{2^{1-\alpha}}{1-2^{\alpha-1}}M_\alpha^o f(x).
\end{align*}
This proves the first assertion.

    Finally, for each $r$,
\[
 |S(x,r)|^{\alpha-1/s}\left(\sum_{S(x,r)}|f|^s\right)^{1/s}
 \le
 \left(\frac{|S(x,r)|}{|B(x,r)|}\right)^{\alpha-1/s}
 |B(x,r)|^{\alpha-1/s}\left(\sum_{B(x,r)}|f|^s\right)^{1/s}.
\]
If $\alpha-1/s\ge0$ the ratio is at most one; if $\alpha-1/s<0$, it is bounded by
$2^{1/s-\alpha}$. Taking the supremum proves
\eqref{eq:Mas-sphere-ball}.

Assume now that $\alpha s<1$. For every $0\le \ell\le r$,
\[
 \left(\sum_{S(x,\ell)}|f|^s\right)^{1/s}
 \le M_{\alpha,s}^o f(x)\,|S(x,\ell)|^{1/s-\alpha}.
\]
Therefore, by \eqref{eq:sphere-size},
\begin{align*}
|B(x,r)|^{\alpha-1/s}
 \left(\sum_{B(x,r)}|f|^s\right)^{1/s}
&\le M_{\alpha,s}^o f(x)\,
 |B(x,r)|^{\alpha-1/s}
 \left(\sum_{\ell=0}^r |S(x,\ell)|^{1-\alpha s}\right)^{1/s}\\
&\lesssim_{\alpha,s} 
 M_{\alpha,s}^o f(x)\,
 k^{r(\alpha-1/s)}k^{r(1/s-\alpha)}
 \lesssim_{\alpha,s}M_{\alpha,s}^o f(x).
\end{align*}
 For the geometric sum to converge, we need $\alpha s <1$.
\end{proof}

\begin{remark}\label{rem:ball-sphere-fractional}
At $\alpha=0$ the first ball--sphere comparison is immediate from the decomposition of a
ball into spheres. In the fractional setting this constant-one argument is no longer
valid because each sphere carries a different power $|S(x,\ell)|^{1-\alpha}$. The
geometric summation in the preceding proof is the first place where the fractional
normalization changes the non-fractional argument.
\end{remark}

We will denote by $\mathbbm{1}\otimes w$ the product measure
\[
 \mathbbm{1}\otimes w(A\times B)
 =|A|\,w(B)
 =\sum_{(x,y)\in A\times B}w(y).
\]

\begin{lemma}\label{1otimesw}
Let $E,F$ be finite subsets of $T$, $s>1$, and $r\in\mathbb N_0$. Then
\begin{equation}\label{eq:incidence-fractional}
\begin{aligned}
\mathbbm{1}\otimes w
\bigl(\{(x,y)\in E\times F:d(x,y)=r\}\bigr)
\le
C_s k^{r\frac{s'}{s'+1}(1-\alpha)}
 w(F)^{\frac1{s'+1}}
 \bigl(M_{\alpha,s}^o w(E)\bigr)^{\frac{s'}{s'+1}},
\end{aligned}
\end{equation}
where $s'=s/(s-1)$ and $C_s$ is independent of $k$.
\end{lemma}

\begin{proof}
Write $T=\bigcup_{j\ge0}T^j$, where $T^j$ is the $j$th generation. 
For each $i,j\geq 0$, we denote  $E_j=E\cap T^j$, $F_i=F\cap T^i$ and 
\[I_{i,j} = \mathbbm{1}\otimes w
\bigl(\{(x,y)\in E_j\times F_i:d(x,y)=r\}\bigr)
\]

If $x\in E_j$, $y\in F_i$ and
$d(x,y)=r$, then
\[
 i=j+r-2m
\]
for some $m\in\{0,\dots,r\}$.
Then $I_{i,j}=0$ if $i\not = j+r-2m$ for any $m\in\{0,\dots,r\}$ and 
\[
\mathbbm{1}\otimes w
\bigl(\{(x,y)\in E\times F:d(x,y)=r\}\bigr) = \sum_{m=0}^r \sum_{\substack{i,j \in \mathbb{N}_0 \\  i=j+r-2m}} I_{i,j}\]

For fixed $y\in T^i$, there are at most $k^m$ vertices $x\in T^j$ at distance $r$,
so
\begin{equation}\label{eq:incidence-first-bound}
 I_{i,j}\le k^m w(F_i).
\end{equation}
On the other hand, for fixed $x\in T^j$ there are at most $k^{r-m}$ vertices of
$T^i$ at distance $r$. By H\"older's inequality and \eqref{eq:sphere-size},
\begin{align*}
I_{i,j}
&=\sum_{x\in E_j} \sum_{\substack{y\in F_i \\ d(x,y)=r}} w(y)
\\&=\sum_{x\in E_j} w(F_i\cap S(x,r))
\\&\leq \sum_{x\in E_j}|F_i\cap S(x,r)|^{\frac1{s'}} w^s(F_i\cap S(x,r))^{\frac1{s}}
\\&= \sum_{x\in E_j}|F_i\cap S(x,r)|^{\frac1{s'}} \frac{k^{r(\frac1{s}-\alpha)}}{k^{r(\frac1{s}-\alpha)}}w^s(F_i\cap S(x,r))^{\frac1{s}}
\\&\leq C_{\alpha,s} \sum_{x\in E_j}|F_i\cap S(x,r)|^{\frac1{s'}} k^{r(\frac1{s}-\alpha)}  |S(x,r)|^{\alpha} \frac{w^s(F_i\cap S(x,r))^{\frac1{s}}}{|S(x,r)|^{\frac1{s}}}
\\&\leq C_{\alpha,s} \sum_{x\in E_j}|F_i\cap S(x,r)|^{\frac1{s'}} k^{r(\frac1{s}-\alpha)} M_{\alpha,s}^o(w)(x)
\\&\leq C_{\alpha,s} k^{\frac{r-m}{s'}} k^{r(\frac1{s}-\alpha)} M_{\alpha,s}^o(w)(E_j)
\end{align*}

We therefore have
\begin{equation}\label{eq:min-incidence}
I_{i,j}
\le C_{\alpha,s}
\min\left\{
 k^{\frac{r-m}{s'}}k^{r(\frac1s-\alpha)}M_{\alpha,s}^o w(E_j),
 k^m w(F_i)
\right\}.
\end{equation}

Set 
\[
 c_j=k^{-j/s'}M_{\alpha,s}^o w(E_j),
 \qquad
 d_i=k^{-i}w(F_i),
 \qquad i,j\geq 0.
\]
Then
\[
 \sum_{j\ge0}k^{j/s'}c_j=M_{\alpha,s}^o w(E),
 \qquad
 \sum_{i\ge0}k^i d_i=w(F).
\]
Since $i=j+r-2m$,
\[
 r-m+j=\frac{r+i+j}{2},
 \qquad
 m+i=\frac{r+i+j}{2}.
\]

Since $i=j+r-2m$, we have
\begin{equation}
\begin{aligned}
\min& \left\{ k^{\frac{r-m}{s'}} k^{r(\frac1{s}-\alpha)} M_{\alpha,s}^o(w)(E_j), k^m w(F_i)\right\}
\\&= \min \left\{ k^{\frac{r-m}{s'}} k^{r(\frac1{s}-\alpha)} k^{\frac{j}{s'}}c_j, k^m k^id_i\right\}
\\&= \min \left\{ k^{\frac{r-m+j}{s'}} k^{r(\frac1{s}-\alpha)} c_j, k^{m+i}d_i \right\}
\end{aligned}
\end{equation}

Since $i=j+r-2m$ then $m=\frac1{2}(j-i+r)$ and 
$$
\frac{r-m+j}{s'}
=\frac1{s'} \left[ r- \frac1{2}(j-i+r) +j\right]
=\frac{ r+ j+i}{2s'} 
$$

Then, for a $\theta$ to be determinate later we have

\begin{align*}
\sum_{m=0}^r & \sum_{\substack{i,j \in \mathbb{N} \cup \{0\} \\ i = j + r - 2m}}  \min \left\{ k^{\frac{r-m+j}{s'}} k^{r(\frac1{s}-\alpha)} c_j, k^{m+i}d_i \right\}
\\&\leq \sum_{m=0}^r \sum_{\substack{i,j \in \mathbb{N} \cup \{0\} \\ i = j + r - 2m}}  \min \left\{ k^{\frac{r}{2s'}} k^{\frac{i+j}{2s'}} k^{r(\frac1{s}-\alpha)} c_j,  k^{\frac{r}{2}} k^{\frac{i+j}{2}} d_i \right\}
\\& \leq \sum_{\substack{i,j\in \mathbb{N}_0 \\ i<j+\theta}} k^{\frac{r}{2s'}} k^{\frac{i+j}{2s'}} k^{r(\frac1{s}-\alpha)} c_j 
+\sum_{\substack{i,j\in \mathbb{N}_0 \\ i\geq j+\theta}} k^{\frac{r}{2}} k^{\frac{i+j}{2}} d_i 
\\&= k^{\frac{r}{2s'}}k^{r(\frac1{s}-\alpha)}  \sum_{\substack{i,j\in \mathbb{N}_0 \\ i<j+\theta}}k^{\frac{i+j}{2s'}}  c_j 
+ k^{\frac{r}{2}}\sum_{\substack{i,j\in \mathbb{N}_0 \\ i\geq j+\theta}} k^{\frac{i+j}{2}} d_i 
\\&= k^{\frac{r}{2s'}}k^{r(\frac1{s}-\alpha)} \sum_{j=0}^{\infty} k^{\frac{j}{2s'}}  c_j  \sum_{\substack{i\in \mathbb{N}_0 \\ i<j+\theta}}k^{\frac{i}{2s'}}  
+ k^{\frac{r}{2}}\sum_{i=0}^{\infty}k^{\frac{i}{2}} d_i  \sum_{\substack{j\in \mathbb{N}_0 \\ j \leq i-\theta}} k^{\frac{j}{2}} 
\\&\leq C_{s}  \left[ k^{\frac{r}{2s'}}k^{r(\frac1{s}-\alpha)} \sum_{j=0}^{\infty} k^{\frac{j}{2s'}}  c_jk^{\frac{j+\theta}{2s'}}  
+  k^{\frac{r}{2}}\sum_{i=0}^{\infty}k^{\frac{i}{2}} d_i  k^{\frac{i-\theta}{2}}\right]
\\&\leq C_{s}  \left[ k^{\frac{r}{2s'}}k^{r(\frac1{s}-\alpha)} k^{\frac{\theta}{2s'}} \sum_{j=0}^{\infty} k^{\frac{j}{s'}}  c_j
+  k^{\frac{r}{2}}k^{-\frac{\theta}{2}} \sum_{i=0}^{\infty}k^{i} d_i  \right]
\\&\leq C_{s}  \left[ k^{\frac{r}{2s'}}k^{r(\frac1{s}-\alpha)} k^{\frac{\theta}{2s'}} M_{\alpha,s}^o(w)(E)
+  k^{\frac{r}{2}}k^{-\frac{\theta}{2}} w(F)  \right]
\end{align*}

where $C_{s}=\frac{2^{\frac1{2s'}}}{2^{\frac1{2s'}}-1}$.

If both terms are nonzero, choose $\theta$ so that they are equal. Equivalently, if
\[
 A=k^{\frac r{2s'}}k^{r(\frac1s-\alpha)}M_{\alpha,s}^o w(E),
 \qquad
 B=k^{\frac r2}w(F),
\]
choose $\theta$ with $A k^{\theta/(2s')}=B k^{-\theta/2}$. Then
\[
 A k^{\theta/(2s')}+B k^{-\theta/2}
 =2A^{\frac{s'}{s'+1}}B^{\frac1{s'+1}}.
\]
The cases $A=0$ or $B=0$ are immediate. Therefore
\begin{align*}
\mathbbm{1}\otimes w
\bigl(\{(x,y)\in E\times F:d(x,y)=r\}\bigr)
\le C_s
 k^{\frac r{s'+1}}
 k^{r(\frac1s-\alpha)\frac{s'}{s'+1}}
 \bigl(M_{\alpha,s}^o w(E)\bigr)^{\frac{s'}{s'+1}}
 w(F)^{\frac1{s'+1}}.
\end{align*}
Since
\[
 \frac1{s'+1}
 +\left(\frac1s-\alpha\right)\frac{s'}{s'+1}
 =\frac{s'}{s'+1}(1-\alpha),
\]
we obtain \eqref{eq:incidence-fractional}.
\end{proof}

\begin{remark}\label{rem:incidence-change}
The exponent $r\frac{s'}{s'+1}(1-\alpha)$ in
\eqref{eq:incidence-fractional} is the first substantial modification of the weighted
incidence estimate from \cite{ORRS}. It comes from the factor
$k^{r(1/s-\alpha)}$ produced when the $L^s$ norm on a sphere is rewritten in terms of
$M_{\alpha,s}^o w$. This is precisely the factor that later interacts with the
renormalization $k^{r\alpha}$ of the function.
\end{remark}

For $0\leq \alpha < 1$ and $r\geq 0$, we define fractional spherical averages 
\[ A_{r,\alpha}^o h(x)
 =\frac1{|S(x,r)|^{1-\alpha}}\sum_{y\in S(x,r)}|h(y)|.
\]
Thus $M_\alpha^o h=\sup_{r\ge0}A_{r,\alpha}^o h$. In the case $\alpha=0$, we denote $A_{r}^o h(x)=A_{r,0}^o h(x)$ the spherical averages.

\begin{lemma}\label{Lemma Avg}
Let $0<\alpha<1$, $s>1$, $r\ge1$, and $\lambda>0$. 
For every finitely supported $f$ and every weight $w$,
\begin{equation}\label{eq:distribution-fractional}
\begin{aligned}
&w\bigl(\{x:A_{r,\alpha}^o f(x)>\lambda\}\bigr)\\
&\quad\le C_{\alpha,s}
\sum_{\substack{n\in\mathbb N_0\\1\le2^n\le k^r}}
\left(\frac{2^n}{k^r}\right)^{1+\frac1{2s'}}
 k^{r(1-\alpha)}
 M_{\alpha,s}^o w
 \left(\left\{x:|f(x)|k^{r\alpha}
 \ge c_\alpha 2^{n-1}\lambda\right\}\right),
\end{aligned}
\end{equation} where $c_{\alpha}$ depends only of $\alpha$ and $C_{\alpha,s}$ depends only of $\alpha$ and $s$. Moreover, $C_{\alpha,s}\to \infty$ when $s\to1$.
\end{lemma}

\begin{proof}
By \eqref{eq:sphere-size},
\[
 A_{r,\alpha}^o f(x)=|S(x,r)|^\alpha A_r^o(|f|)(x).
\]
Hence
\[
 A_{r,\alpha}^o f(x)>\lambda
 \quad\Longrightarrow\quad
 A_r^o g(x)>1,
 \qquad
 g=\frac{2^\alpha k^{r\alpha}}{\lambda}|f|.
\]
We estimate the latter level set.

Let $N_r$ be the largest integer such that $2^{N_r}\le k^r$, and put
\[
 E_n=\{x:2^{n-1}\le g(x)<2^n\},
 \qquad 0\le n\le N_r.
\]
Then
\begin{equation}\label{eq:decomp-g}
 g\le\frac12+\sum_{n=0}^{N_r}2^n\chi_{E_n}
 +g\chi_{\{g\ge k^r/2\}}.
\end{equation}
Consequently,
\begin{equation}\label{eq:split-level-g}
\{A_r^o g>1\}
\subset
\left\{\sum_{n=0}^{N_r}2^nA_r^o(\chi_{E_n})\ge\frac12\right\}
\cup \left\{A_r^o\bigl(g\chi_{\{g\ge k^r/2\}}\bigr)\ne0\right\},
\end{equation}

Let $0<\beta<1$ to be chosen later. If
\[
 \sum_{n=0}^{N_r}2^nA_r^o(\chi_{E_n})(x)\ge\frac12,
\]
then for some $0\le n\le N_r$,
\begin{equation}\label{eq:Fn-threshold}
 A_r^o(\chi_{E_n})(x)
 \ge\frac{2^\beta-1}{2^{n+2}}
 \left(\frac{2^n}{k^r}\right)^\beta.
\end{equation}
Indeed, if \eqref{eq:Fn-threshold} failed for all $n$, then
\begin{align*}
\sum_{n=0}^{N_r}2^nA_r^o(\chi_{E_n})(x)
&<\frac{2^\beta-1}{4k^{r\beta}}
 \sum_{n=0}^{N_r}2^{\beta n}\\
&<\frac{2^{\beta(N_r+1)}}{4k^{r\beta}}
\le\frac{2^\beta}{4}<\frac12,
\end{align*}
a contradiction.

Define
\[
 F_n=\left\{x:A_r^o(\chi_{E_n})(x)\ge\tau_n\right\},
 \qquad
 \tau_n=\frac{2^\beta-1}{2^{n+2}}
 \left(\frac{2^n}{k^r}\right)^\beta.
\]
Since $f$ is finitely supported, $E_n$ and $F_n$ are finite. For $y\in F_n$,
\[
 |E_n\cap S(y,r)|\ge\tau_n|S(y,r)|\ge\tau_n k^r.
\]
Therefore
\begin{equation}\label{eq:pair-lower}
\mathbbm{1}\otimes w
\bigl(\{(x,y)\in E_n\times F_n:d(x,y)=r\}\bigr)
\ge\tau_n k^r w(F_n).
\end{equation}
By Lemma~\ref{1otimesw},
\begin{equation}\label{eq:pair-upper}
\begin{aligned}
&\mathbbm{1}\otimes w
\bigl(\{(x,y)\in E_n\times F_n:d(x,y)=r\}\bigr)\\
&\quad\le C_s k^{r\frac{s'}{s'+1}(1-\alpha)}
 w(F_n)^{\frac1{s'+1}}
 \bigl(M_{\alpha,s}^o w(E_n)\bigr)^{\frac{s'}{s'+1}}.
\end{aligned}
\end{equation}
Combining \eqref{eq:pair-lower} and \eqref{eq:pair-upper}, and raising to the power
$(s'+1)/s'$, gives
\begin{equation}\label{eq:Fn-general-beta}
 w(F_n)
 \le C_{s,\beta}
 \left(\frac{2^n}{k^r}\right)^{(1-\beta)\frac{s'+1}{s'}}
 k^{r(1-\alpha)}M_{\alpha,s}^o w(E_n).
\end{equation}
Choose
\[
 \beta=\frac1{2(s'+1)}.
\]
Then
\[
 (1-\beta)\frac{s'+1}{s'}=1+\frac1{2s'},
\]
so
\begin{align}\label{eq:Fn-final}
 w(F_n)
 &\le C_s
 \left(\frac{2^n}{k^r}\right)^{1+\frac1{2s'}}
 k^{r(1-\alpha)}M_{\alpha,s}^o w(E_n) \nonumber
\\&\le C_s
 \left(\frac{2^n}{k^r}\right)^{1+\frac1{2s'}}
 k^{r(1-\alpha)}M_{\alpha,s}^o w(\{g\ge2^{n-1}\})
\end{align}

It remains to control of $\left\{A_r^o\bigl(g\chi_{\{g\ge k^r/2\}}\bigr)\ne0\right\}$. 
From its definition,
\[
 \left\{A_r^o\bigl(g\chi_{\{g\ge k^r/2\}}\bigr)\ne0\right\} \subset \bigcup_{y\in\{g\ge k^r/2\}}S(y,r).
\]
By H\"older's inequality,
\begin{align*}
 w(S(y,r))
&\le |S(y,r)|^{1/s'}
 \left(\sum_{z\in S(y,r)}w(z)^s\right)^{1/s}\\
&=|S(y,r)|^{1-\alpha}
 \left[|S(y,r)|^\alpha
 \left(\frac1{|S(y,r)|}\sum_{z\in S(y,r)}w(z)^s\right)^{1/s}\right]\\
&\le C_\alpha k^{r(1-\alpha)}M_{\alpha,s}^o w(y).
\end{align*}
Thus
\begin{equation}\label{eq:tail-estimate}
 w\left(\left\{A_r^o\bigl(g\chi_{\{g\ge k^r/2\}}\bigr)\ne0\right\}\right)
 \le C_\alpha k^{r(1-\alpha)}M_{\alpha,s}^o w\left(\left\{g\ge k^r/2\right\}\right).
\end{equation}
Now $2^{N_r}\le k^r<2^{N_r+1}$, hence
\[
 \{g\ge k^r/2\}\subset\{g\ge2^{N_r-1}\},
 \qquad
 \frac{2^{N_r}}{k^r}>\frac12.
\]
Therefore \eqref{eq:tail-estimate} is dominated, up to a constant, by the term
$n=N_r$ in the sum obtained from \eqref{eq:Fn-final}. This is the point where the cutoff
$2^n\le k^r$ is essential.

Combining this observation with \eqref{eq:split-level-g},
\[
 w(\{A_r^o g>1\})
 \le C_{\alpha,s}\sum_{n=0}^{N_r}
 \left(\frac{2^n}{k^r}\right)^{1+\frac1{2s'}}
 k^{r(1-\alpha)}M_{\alpha,s}^o w(\{g\ge2^{n-1}\}).
\]
Since $g=2^\alpha k^{r\alpha}|f|/\lambda$, this is
\eqref{eq:distribution-fractional}, for instance with $c_\alpha=2^{-\alpha}$.
The explicit constants in \eqref{eq:Fn-general-beta}, together with
$\beta\downarrow0$ as $s\downarrow1$, show that the resulting constant may be chosen
to satisfy $C_{\alpha,s}\to\infty$ as $s\downarrow1$.
\end{proof}

\begin{remark}\label{rem:cutoff-novelty}
Following the account in \cite{ORRS}, the first attempt is to cut the dyadic levels in $k^{r\alpha}$ instead of $k^r$. That cutoff is too short: its last dyadic level is not
comparable with the geometric size $k^r$, and the large-value tail cannot be absorbed.
The correct procedure is instead to first normalize $f$ by $k^{r\alpha}$ and then retain
the geometric cutoff $k^r$ from the underlying sphere. This distinction is invisible
when $\alpha=0$ and is one of the key structural points of the fractional proof.
\end{remark}

Now, we proceed to the proof of the main result, Theorem \ref{teo fS con s}.

\begin{proof}[Proof of Theorem \ref{teo fS con s}]
We first assume that $f$ is finitely supported. By Proposition~\ref{Malpho}, it is enough 
to prove the corresponding estimate for $M_\alpha^o$, with $M_{\alpha,s}^o w$ on the
right-hand side.
Observe that 
\[
w(\{M_\alpha^o f>\lambda\}) \leq \sum_{r=0}^{\infty} w(\{A_{r,\alpha}^o f>\lambda\}).
\]

The radius $r=0$ contributes
\[
 w(\{|f|>\lambda\})
 \le\frac1\lambda\sum_{x\in T}|f(x)|w(x)
 \le\frac1\lambda\sum_{x\in T}|f(x)|M_{\alpha,s}^o w(x).
\]

For $r\ge1$, Lemma~\ref{Lemma Avg} gives,
\begin{align*}
&w(\{M_\alpha^o f>\lambda\})\\
&\quad\le C_{\alpha,s}
\sum_{r\ge1}\sum_{\substack{n\ge0\\2^n\le k^r}}
\left(\frac{2^n}{k^r}\right)^{1+\frac1{2s'}}
 k^{r(1-\alpha)}
 M_{\alpha,s}^o w
 \bigl(\{|f|k^{r\alpha}\ge c_\alpha2^{n-1}\lambda\}\bigr)
 +\frac{C_\alpha}{\lambda}\sum|f|M_{\alpha,s}^o w.
\end{align*}
Interchanging the sums, fix $x\in T$ and set
\[
 A=C_\alpha\frac{|f(x)|}{\lambda},
 \qquad
 a=1+\frac1{2s'}>1.
\]
The contribution associated with $x$ is bounded by
\begin{align*}
\Phi(A)
&:=\sum_{r\ge1}k^{r(1-\alpha)}
 \sum_{\substack{n\ge0\\2^n\le k^r\\2^n\le A k^{r\alpha}}}
 \left(\frac{2^n}{k^r}\right)^a\\
&\lesssim_s
 \sum_{r\ge1}k^{r(1-\alpha-a)}
 \min\{k^r,A k^{r\alpha}\}^{a}\\
&=\sum_{r\ge1}k^{r(1-\alpha)}
 \min\{1,A k^{-r(1-\alpha)}\}^{a}.
\end{align*}
We claim that
\begin{equation}\label{eq:fractional-summation}
 \Phi(A)\lesssim_{\alpha,s}A,
 \qquad A\ge0,
\end{equation}
with a constant independent of $k\ge2$. If $0\le A\le1$, then
\[
 \Phi(A)
 \le A^a\sum_{r\ge1}k^{-r(1-\alpha)(a-1)}
 \lesssim_{\alpha,s}A^a
 \le A.
\]
If $A>1$, choose $R\in\mathbb N_0$ so that
\[
 k^{R(1-\alpha)}\le A<k^{(R+1)(1-\alpha)}.
\]
Then
\begin{align*}
\Phi(A)
&\le\sum_{r\le R}k^{r(1-\alpha)}
 +A^a\sum_{r>R}k^{-r(1-\alpha)(a-1)}
 \lesssim_{\alpha,s}A.
\end{align*}
Since $k\ge2$, all geometric-series constants are bounded uniformly in $k$. This proves
\eqref{eq:fractional-summation}.

It follows that
\[
 w(\{M_\alpha^o f>\lambda\})
 \le C_{\alpha,s}\frac1\lambda
 \sum_{x\in T}|f(x)|M_{\alpha,s}^o w(x).
\]
Using Proposition~\ref{Malpho} once more, namely
$M_\alpha\lesssim_\alpha M_\alpha^o$ and
$M_{\alpha,s}^o w\le2^{(1/s-\alpha)_+}M_{\alpha,s}w$, gives
\eqref{eq:FS-main}.

For a general $f$, apply the estimate to an increasing sequence of finitely supported
truncations of $|f|$ and use monotone convergence on both sides.
\end{proof}

\begin{remark}\label{rem:new-summation}
The proof of Theorem~\ref{teo fS con s} contains a second genuinely fractional
feature. After the dyadic sum, the radius sum does not reduce to the same geometric
series as in the case $\alpha=0$. Instead, the two scales $k^{r\alpha}$ and
$k^{r(1-\alpha)}$ combine into the kernel
\[
 k^{r(1-\alpha)}\min\{1,A k^{-r(1-\alpha)}\}^{1+\frac1{2s'}},
\]
and the strict inequality $1+\frac1{2s'}>1$ is exactly what makes the tail summable. This gives a
fractional analogue of the decay mechanism emphasized in \cite{ORRS}, but with a new
radius--level interaction.
\end{remark}

 As an immediate result of Theorem \ref{teo fS con s} we obtain the following sufficient condition for the weighted weak-type $(1,1)$ of the fractional maximal.

\begin{corol}\label{cor:weighted-weak}
Let $k\ge2$ and $0<\alpha<1$. Assume that there exist $s>1$ and $C_w>0$ such that
\[
 M_{\alpha,s}w(x)\le C_w w(x),\qquad x\in T.
\]
Then
\[
 w\bigl(\{x\in T:M_\alpha f(x)>\lambda\}\bigr)
 \lesssim_{\alpha,s,C_w}
 \frac1\lambda\sum_{x\in T}|f(x)|w(x),
\]
with a constant independent of $k$.
\end{corol}

\section{Radial weights, sharpness, and necessary conditions}

Let $T$ be the infinite rooted $k$-ary tree and write
\[
T=\bigcup_{j=0}^{\infty}T^j,
\]
where $T^0$ consists only of the root and $T^j$ is the generation at depth $j$.
We consider radial weights of the form
\[
w(x)=\sum_{j=0}^{\infty}c_j\chi_{T^j}(x),\qquad c_j\geq0.
\]

\subsection{Radial weights and the subcritical range}

A natural first problem is to identify explicit weights satisfying the fractional
$A_1$-type condition
\[
 M_{\alpha}w\lesssim w.
\]

One would expect the natural response to be that for each $0\leq \alpha < 1$ find a power weight depending on $\alpha$. When $\alpha=0$ should be coincide with the one provided in \cite{ORRS}. However, the
admissible interval of radial powers shrinks as $\alpha$ increases and disappears at the
critical value $\alpha=1/2$.

For $\beta\in\mathbb R$, define
\[
 w_{\beta}(x)=\sum_{j=0}^{\infty}k^{j\beta}\chi_{T^j}(x).
\]
We first record the geometric computation that will be used later.

\begin{lemma}\label{lem:radial-sphere-decomp}
Let $x\in T^j$, $r\ge0$, and $L=\min\{j,r\}$. Then, uniformly in
$k\ge2$, $j$ and $r$,
\begin{equation}\label{eq:radial-frac}
 \frac1{|S(x,r)|^{1-\alpha}}
 \sum_{y\in S(x,r)}w_\beta(y)
 \simeq_{\alpha} w_{\beta}(x)\,
 k^{r(\alpha+\beta)}
 \sum_{m=0}^{L}k^{-m(1+2\beta)}.
\end{equation}
More generally, for $1\le s<\infty$,
\begin{equation}\label{eq:radial-Mas}
|S(x,r)|^{\alpha}
 \left( \frac1{|S(x,r)|}\sum_{y\in S(x,r)}w_\beta^s(y) \right)^{1/s}
 \simeq_{\alpha,s} w_\beta(x)\, 
 k^{r(\alpha+\beta)}
 \left(\sum_{m=0}^{L}k^{-m(1+2s\beta)}\right)^{1/s},
\end{equation}
\end{lemma}

\begin{proof}
Let $x\in T^j$, if $y\in S(x,r) \cap T^i$ then $i=j+r-2m$ for some $0\le m\le L$. Then splitting 
\[S(x,r) = \bigcup_{m=0}^L \{ y \in S(x,r) : y\in T^i, i=j+r-2m \}.\]
The number of vertices in the corresponding
class is comparable, uniformly in $k\ge2$, to $k^{r-m}$. Indeed, for $m=0$ there are
exactly $k^r$ descendants; for an interior class the exact number is
$(k-1)k^{r-m-1}$; and the possible terminal ancestor contributes one point. Hence
\begin{align*}
 \sum_{y\in S(x,r)}w_\beta(y)
 &\simeq
 \sum_{m=0}^{L}k^{r-m}k^{\beta(j+r-2m)}\\
 &=k^{j\beta}k^{r(1+\beta)}
 \sum_{m=0}^{L}k^{-m(1+2\beta)}.
\end{align*}
Since $|S(x,r)|\simeq k^r$, this gives \eqref{eq:radial-frac}. Replacing
$w_\beta$ by $w_\beta^s$ in the unnormalized sum, taking the $s$-th root and
multiplying by $|S(x,r)|^{\alpha-1/s}$ gives \eqref{eq:radial-Mas}.
\end{proof}

We can now characterize the radial powers completely.

\begin{prop}\label{prop:radial-weights}
Let $0<\alpha<1$. If $0<\alpha<1/2$, then
\[
 M_{\alpha}w_{\beta}\simeq_{\alpha,\beta}w_{\beta}
\]
uniformly on $T$, with constants independent of $k\ge2$, if and only if
\[\alpha-1\le\beta\le-\alpha.
\]
If $1/2\le\alpha<1$, there is no $\beta\in\mathbb R$ for which
$M_{\alpha}w_{\beta}\lesssim w_{\beta}$ holds.
\end{prop}

\begin{proof}
Since the radius $r=0$ gives $M_\alpha w_\beta\ge w_\beta$, it remains to determine
when the reverse inequality holds. By Proposition~\ref{Malpho}, it is enough to work
with $M_\alpha^o$.
By Lemma \ref{lem:radial-sphere-decomp}, for $x\in T^j$ we have 
\begin{equation*}
A_{r,\alpha}^{\circ}(w_{\beta})(x)= \frac1{|S(x,r)|^{1-\alpha}}
 \sum_{y\in S(x,r)}w_\beta(y)
 \simeq_{\alpha} w_{\beta}(x)\,
 k^{r(\alpha+\beta)}
 \sum_{m=0}^{L}k^{-m(1+2\beta)}.
\end{equation*}

First assume $j\ge r$, so that $L=r$ in Lemma  \ref{lem:radial-sphere-decomp}. If
$\beta<-1/2$, then
\[
 \sum_{m=0}^{r}k^{-m(1+2\beta)}
 \simeq_\beta k^{r(-1-2\beta)},
\]
and therefore 
\[
\frac{A_{r,\alpha}^o(w_{\beta})(x)}{w_{\beta}(x)}  \simeq k^{r(\alpha-\beta-1)}.
\]

Uniform boundedness forces $\beta\ge\alpha-1$. If $\beta=-1/2$, 
\[ \frac{A_{r,\alpha}^o(w_{\beta})(x)}{w_{\beta}(x)} \simeq 
 (r+1)k^{r(\alpha-1/2)},
\]
which is uniformly bounded exactly when $\alpha<1/2$. Finally, if
$\beta>-1/2$, the geometric sum is uniformly bounded and 
\[\frac{A_{r,\alpha}^o(w_{\beta})(x)}{w_{\beta}(x)} \simeq k^{r(\alpha+\beta)}, \]
 so necessarily $\beta\le-\alpha$.

These conditions are also sufficient. In the first case, the conditions
$\beta<-1/2$ and $\beta\ge\alpha-1$ imply $\alpha<1/2$ and
$\alpha+\beta<0$. If $r\le j$, the preceding estimate gives
\[
 k^{r(\alpha-\beta-1)}\le1.
\]
If $r>j$, then $L=j$ and
\[
 k^{r(\alpha+\beta)}k^{j(-1-2\beta)}
 \le k^{j(\alpha-\beta-1)}\le1.
\]
At $\beta=-1/2$, the strict inequality $\alpha<1/2$ gives
\[
 (L+1)k^{r(\alpha-1/2)}\lesssim_\alpha1
\]
for both $r\le j$ and $r>j$. In the last case the geometric sum is bounded and
$\alpha+\beta\le0$. Thus $M_\alpha^o w_\beta\lesssim w_\beta$, and the ball
estimate follows from Proposition~\ref{Malpho}.

Combining the three regimes gives exactly $[\alpha-1,-\alpha]$ when
$0<\alpha<1/2$. At $\alpha=1/2$ the only possible meeting point would be
$\beta=-1/2$, but the factor $r+1$ prevents boundedness; for $\alpha>1/2$ the two
necessary inequalities are incompatible.
\end{proof}

The behavior of $L^s$ fractional maximal operator $M_{\alpha,s}$, $s>1$, we have

\begin{prop}\label{prop:radial-Mas}
Let $0<\alpha<1$ and $s>1$. Then
\[
 M_{\alpha,s}w_\beta\simeq_{\alpha,s,\beta}w_\beta
\]
uniformly on $T$, with constants independent of $k\ge2$, if and only if
\[
 s< \frac1{2\alpha},\qquad
 \alpha-\frac1s\le\beta\le-\alpha.
\]
In particular, nontrivial radial weights satisfying
$M_{\alpha,s}w\lesssim w$ exist if and only if
$0<\alpha<1/2$ and $1<s<1/(2\alpha)$.
\end{prop}

\begin{proof} The proof is analogous to the case $s=1$ proved above, using the estimate \eqref{eq:radial-Mas} instead of \eqref{eq:radial-frac}. We provide the complete proof to show the differences. 

Necessity follows already from the spherical operator, because
$M_{\alpha,s}^o\lesssim M_{\alpha,s}$ by \eqref{eq:Mas-sphere-ball}. Take
$j\ge r$ in \eqref{eq:radial-Mas}. If $\beta<-1/(2s)$, then
\[
 \left(\sum_{m=0}^{r}k^{-m(1+2s\beta)}\right)^{1/s}
 \simeq_{s,\beta}k^{r(-1/s-2\beta)},
\]
so the quotient is comparable to $k^{r(\alpha-\beta-1/s)}$. Hence
$\beta\ge\alpha-1/s$. If $\beta=-1/(2s)$, the quotient is
\[
 (r+1)^{1/s}k^{r(\alpha-1/(2s))},
\]
and is bounded only if $\alpha<1/(2s)$. If $\beta>-1/(2s)$, the geometric
sum is bounded and we must have $\beta\le-\alpha$. Thus the necessary interval is
nonempty precisely when $\alpha\leq 1/(2s)$, or equivalently
$s\leq 1/(2\alpha)$. 
At equaly, its only point is $\beta=-\frac1{2s}$, which  is excluded by the proceding computation. Hence the admissible range is $\alpha<\frac1{2s}$

Conversely, assume these conditions. In particular $\alpha s<1/2<1$, so by
\eqref{eq:Mas-ball-sphere} it suffices to bound $M_{\alpha,s}^o w_\beta$.
When $\beta<-1/(2s)$ and $r>j$, formula \eqref{eq:radial-Mas} gives
\[
 k^{r(\alpha+\beta)}k^{j(-1/s-2\beta)}.
\]
Here $\alpha+\beta<0$, and hence this is bounded by
\[
 k^{j(\alpha-\beta-1/s)}\le1.
\]
For $r\le j$ the same conclusion follows from
$\alpha-\beta-1/s\le0$. At the critical value $\beta=-1/(2s)$, the factor
$(L+1)^{1/s}$ is dominated by the exponential decay coming from
$\alpha-1/(2s)<0$. Finally, if $\beta>-1/(2s)$, the geometric sum is uniformly
bounded and $\alpha+\beta\le0$. This proves the spherical estimate and hence, by
\eqref{eq:Mas-ball-sphere}, the ball estimate. The reverse inequality follows from
$r=0$.
\end{proof}

\begin{remark}\label{rem:radial-consistency}
These formulas recover the setting $\alpha=0$ in Proposition~\ref{prop:radial-weights} gives the interval
$[-1,0]$ from \cite[Theorem~3.1]{ORRS}. Similarly,
Proposition~\ref{prop:radial-Mas} becomes
$ -\frac1s\le\beta\le0.$

Thus, every $\beta\in(-1,0]$ admits some $s>1$ sufficiently close to $1$, as the
extreme power $\beta=-1$ admits no such $s>1$, in agreement with
\cite[Remarks~3.2--3.3]{ORRS}.

For $\alpha>0$ the behavior is genuinely different. The interval
$[\alpha-1,-\alpha]$ has length $1-2\alpha$ and collapses at
$\alpha=1/2$. For $M_{\alpha,s}$ the admissible interval has length
$1/s-2\alpha$, which explains the sharp restriction $s<1/(2\alpha)$. These shrinking
intervals play a crucial role in the behavior changes at $\alpha=\frac1{2}$.
\end{remark}

\begin{remark}\label{rem:nontrivial-weak}
If $0<\alpha<1/2$, choose $1<s<1/(2\alpha)$ and
$\beta\in[\alpha-1/s,-\alpha]$. Proposition~\ref{prop:radial-Mas} and
Corollary~\ref{cor:weighted-weak} then give explicit nonconstant weights for which
\[
 w_\beta(\{M_\alpha f>\lambda\})
 \lesssim_{\alpha,s,\beta}
 \frac1\lambda\sum_{x\in T}|f(x)|w_\beta(x).
\]
Thus the positive Fefferman--Stein estimate is non-vacuous in a genuinely subcritical
weighted regime.
\end{remark}

\subsection{Failure of the classical fractional Fefferman--Stein estimate}

We now show that the estimate obtained above with $M_{\alpha,s}w$, $s>1$, cannot in
general be extended to $s=1$. More precisely, we construct a weight $w$ and a
sequence $\{f_j\}$ for which
\[
w\bigl(\{x\in T:M_{\alpha}f_j(x)>1\}\bigr)\longrightarrow\infty,
\]
while
\[
\sup_j\sum_{x\in T}|f_j(x)|M_{\alpha}w(x)<\infty.
\]

The main idea is to construct a subtree whose growth is of order $k^{(1-\alpha) r}$,
which exactly compensates the factor $|B(x,r)|^{1-\alpha}$ appearing in the
fractional maximal operator.

\subsubsection*{A subtree of fractional growth}

For $\ell\geq0$, define
\[
\varepsilon_{\ell}
=\lfloor(\ell+1)(1-\alpha)\rfloor-\lfloor\ell(1-\alpha)\rfloor.
\]
Since $0<\alpha<1$, we have $\varepsilon_{\ell}\in\{0,1\}$.
We construct a rooted subtree $\mathcal S\subset T$, with the same root of $T$, recursively as follows.
Every vertex of $\mathcal S$ in generation $\ell$ keeps $k^{\varepsilon_{\ell}}$
children. Thus, if $\varepsilon_{\ell}=1$, all its $k$ children are kept, whereas
if $\varepsilon_{\ell}=0$, only one child is kept.

Let
\[
\mathcal S_j=\mathcal S\cap T^j,
\qquad
N_j=|\mathcal S_j|.
\]
Then
\begin{equation}\label{eq:Nj}
N_j
=\prod_{\ell=0}^{j-1}k^{\varepsilon_{\ell}}
=k^{\lfloor j(1-\alpha)\rfloor}.
\end{equation}
In particular,
\begin{equation}\label{eq:Nj-growth}
k^{-1}k^{(1-\alpha)j}\leq N_j\leq k^{(1-\alpha) j}.
\end{equation}
Moreover, if $x\in\mathcal S_i$, then the number of descendants of $x$ belonging
to $\mathcal S_{i+r}$ is exactly
\begin{equation}\label{eq:descendants}
D(i,r)=\frac{N_{i+r}}{N_i}
=k^{\lfloor(i+r)(1-\alpha)\rfloor-\lfloor i(1-\alpha)\rfloor},
\end{equation}
and therefore
\begin{equation}\label{eq:descendants-growth}
k^{-1}k^{(1-\alpha) r}\leq D(i,r)\leq k k^{(1-\alpha) r}.
\end{equation}

We define the weight
\begin{equation}\label{eq:subtree-weight}
w(x)=\sum_{j=0}^{\infty}\frac1{N_j}\chi_{\mathcal S_j}(x).
\end{equation}
Thus $w=0$ outside $\mathcal S$ and
\begin{equation}\label{eq:level-weight}
w(\mathcal S_j)=1
\qquad\text{for every }j\geq0.
\end{equation}

The next lemma is the key estimate for the right-hand side of the counterexample.

\begin{lemma}\label{lemma:subtree-weight}
For every $x\in\mathcal S$,
\[
M_{\alpha}^{o}w(x)\lesssim_{\alpha,k}w(x)
\qquad\text{and}\qquad
M_{\alpha}w(x)\lesssim_{\alpha,k}w(x).
\]
Moreover,
\[
M_{\alpha}^{o}w(x)\simeq_{\alpha,k}w(x),
\qquad
M_{\alpha}w(x)\simeq_{\alpha,k}w(x),
\qquad x\in\mathcal S.
\]
\end{lemma}

\begin{proof}
Fix $x\in\mathcal S_i$ and $r\geq0$. For
$0\leq m\leq\min\{r,i\}$, let $E_m$ be the set of vertices
$y\in\mathcal S\cap S(x,r)$ such that the lowest common ancestor of $x$ and $y$
is the ancestor of $x$ at generation $i-m$.

If $y\in E_m$, then $y$ belongs to generation $i+r-2m$. Let $z$ be the ancestor
of $x$ at generation $i-m$. The total number of descendants of $z$ in
$\mathcal S_{i+r-2m}$ is
\[
\frac{N_{i+r-2m}}{N_{i-m}}.
\]
Hence
\[
|E_m|\leq\frac{N_{i+r-2m}}{N_{i-m}}.
\]
Since every vertex of $\mathcal S_{i+r-2m}$ has weight
$N_{i+r-2m}^{-1}$, it follows that
\[
\sum_{y\in E_m}w(y)\leq\frac1{N_{i-m}}.
\]
Therefore,
\begin{equation}\label{eq:sphere-w-sum}
\sum_{y\in S(x,r)}w(y)
\leq\sum_{m=0}^{\min\{r,i\}}\frac1{N_{i-m}}.
\end{equation}

Using \eqref{eq:Nj},
\[
\frac{N_i}{N_{i-m}}
=k^{\lfloor(1-\alpha) i\rfloor-\lfloor(1-\alpha)(i-m)\rfloor}
\leq k k^{(1-\alpha) m},
\]
and consequently
\[
\frac1{N_{i-m}}\leq k k^{(1-\alpha) m}\frac1{N_i}.
\]
It follows from \eqref{eq:sphere-w-sum} that
\begin{equation}\label{eq:sphere-w-final}
\sum_{y\in S(x,r)}w(y)
\lesssim_{\alpha,k}\frac{k^{(1-\alpha) r}}{N_i}.
\end{equation}

Since $x$ has exactly $k^r$ descendants at distance $r$,
$|S(x,r)|\geq k^r$. Hence
\begin{align*}
\frac1{|S(x,r)|^{1-\alpha}}\sum_{y\in S(x,r)}w(y)
&\lesssim_{\alpha,k}
\frac1{k^{(1-\alpha) r}}\frac{k^{(1-\alpha) r}}{N_i}
= \frac1{N_i}
=w(x).
\end{align*}
Taking the supremum in $r$ gives
\[
M_{\alpha}^{o}w(x)\lesssim_{\alpha,k}w(x).
\]

For the maximal function over balls, using \eqref{eq:sphere-w-final},
\begin{align*}
\sum_{y\in B(x,r)}w(y)
&=\sum_{\ell=0}^{r}\sum_{y\in S(x,\ell)}w(y)
\lesssim_{\alpha,k}\frac1{N_i}\sum_{\ell=0}^{r}k^{(1-\alpha)\ell}
\lesssim_{\alpha,k}\frac{k^{(1-\alpha) r}}{N_i}.
\end{align*}
Since $|B(x,r)|\geq k^r$,
\[
\frac1{|B(x,r)|^{1-\alpha}}\sum_{y\in B(x,r)}w(y)
\lesssim_{\alpha,k}\frac1{N_i}=w(x).
\]
Thus
\[
M_{\alpha}w(x)\lesssim_{\alpha,k}w(x).
\]
Finally, taking $r=0$ in either maximal operator yields the reverse inequalities
$M_{\alpha}^{o}w(x)\geq w(x)$ and $M_{\alpha}w(x)\geq w(x)$.
\end{proof}

We now prove the counterexample. A function supported on a single generation
of the subtree has uniformly bounded pairing with $M_\alpha w$, whereas its
maximal function exceeds a fixed level on every preceding generation.
Since each of these generations has unit weight, the distribution function
grows without bound. This construction works for every $0<\alpha<1$.

\begin{teo}\label{teo:failure-classical-FS}
Let $k\geq2$ and $0<\alpha<1$. There exist a weight $w\geq0$ on $T$ and a
sequence of finitely supported functions $\{f_j\}_{j\geq1}$ such that
\begin{equation}\label{eq:RHS-bounded}
\sup_{j\geq1}\sum_{x\in T}|f_j(x)|M_{\alpha}w(x)<\infty,
\end{equation}
whereas
\begin{equation}\label{eq:LHS-unbounded}
w\bigl(\{x\in T:M_{\alpha}f_j(x)>1\}\bigr)\longrightarrow\infty
\qquad\text{as }j\to\infty.
\end{equation}
Consequently, there is no constant $C>0$ such that
\[
w\bigl(\{x\in T:M_{\alpha}f(x)>\lambda\}\bigr)
\leq\frac{C}{\lambda}\sum_{x\in T}|f(x)|M_{\alpha}w(x)
\]
holds for every finitely supported $f$ and every $\lambda>0$.
\end{teo}

\begin{proof}
Let $w$ be the weight defined in \eqref{eq:subtree-weight}. For every $j\geq1$ set
\[
f_j=C_0\chi_{\mathcal S_j},
\]
where $C_0>4^{1-\alpha}k$ is fixed.

We first prove that
\begin{equation}\label{eq:inclusion-counterexample}
\bigcup_{i=0}^{j}\mathcal S_i
\subset
\{x\in T:M_{\alpha}f_j(x)>1\}.
\end{equation}
Let $x\in\mathcal S_i$, $0\leq i\leq j$, and denote $r=j-i$.
By construction, the number of descendants of $x$ in $\mathcal S_j$ is exactly
$D(i,r)=N_j/N_i$. By \eqref{eq:descendants-growth},
\[
D(i,r)\geq k^{-1}k^{(1-\alpha) r}.
\]
All these descendants belong to $B(x,r)$, and therefore
\begin{equation}\label{eq:fj-ball-sum}
\sum_{y\in B(x,r)}f_j(y)
\geq C_0D(i,r)
\geq C_0k^{-1}k^{(1-\alpha) r}.
\end{equation}

Recall that $|S(x,\ell)|\leq 2 k^{\ell}$ and $|B(x,\ell)|\leq 4 k^{\ell}$.
Putting $\delta=1-\alpha$ and combining this with \eqref{eq:fj-ball-sum},
\begin{align*}
M_{\alpha}f_j(x)
&\geq\frac1{|B(x,r)|^{\delta}}
\sum_{y\in B(x,r)}f_j(y)
\geq
\frac{C_0k^{-1}k^{\delta r}}{4^{\delta}k^{\delta r}}
=\frac{C_0}{4^{1-\alpha}k}>1.
\end{align*}
This proves \eqref{eq:inclusion-counterexample}.

Using $w(\mathcal S_i)=1$ for all $i\geq 0$, we obtain
\begin{align*}
w\bigl(\{M_{\alpha}f_j>1\}\bigr)
\geq w\left(\bigcup_{i=0}^{j}\mathcal S_i\right)
=\sum_{i=0}^{j}w(\mathcal S_i)=j+1.
\end{align*}
Hence, the left-hand side tends to infinity.

On the other hand, since $f_j$ is supported on $\mathcal S_j$, Lemma \ref{lemma:subtree-weight} yields
\begin{align*}
\sum_{x\in T}|f_j(x)|M_{\alpha}w(x)
&=C_0\sum_{x\in\mathcal S_j}M_{\alpha}w(x)
\lesssim_{\alpha,k}C_0\sum_{x\in\mathcal S_j}w(x)
=C_0w(\mathcal S_j)=C_0.
\end{align*}
This bound is uniform in $j$, and the conclusion follows.
\end{proof}

\begin{remark}
The essential feature of the construction is that, for $x\in\mathcal S_i$, the
number of descendants of $x$ in $\mathcal S_{i+r}$ satisfies
\[
D(i,r)\simeq_k k^{(1-\alpha)r}.
\]
This growth exactly compensates the denominator
$|B(x,r)|^{1-\alpha}\simeq k^{(1-\alpha)r}$ in the fractional maximal operator.
This is the fractional counterpart of the mechanism used in the non-fractional
counterexample of \cite{ORRS}.
\end{remark}

\subsection{A global counterexample and failure under finite iterations}

We now focus on Theorem \ref{Teo 1.2} in the subcritical range $0<\alpha<1/2$.
For $n\geq1$, let $M_\alpha^{(n)}$ denote the $n$-fold iteration of $M_\alpha$.

\begin{proof}[Proof of Theorem \ref{Teo 1.2}]

Consider the weight 
\[
w_*(x)=\sum_{j=0}^{\infty} k^{(\alpha-1)j} \chi_{T^{j}}(x).
\]

The first assertion follows from Proposition~\ref{prop:radial-weights}, since
$w_*=w_{\alpha-1}$ and $\alpha-1\in[\alpha-1,-\alpha]$ when
$0<\alpha<1/2$. Hence, for each fixed $n$,
\begin{equation}\label{eq:iterate-wstar}
M_\alpha^{(n)}w_*\lesssim_{\alpha,k,n}w_*.
\end{equation}

Let $\mathcal S$ be the subtree constructed above and set
\[
f_j=C_0\chi_{\mathcal S_j},
\]
with the same fixed constant $C_0>4^{1-\alpha}k$. The proof of
Theorem~\ref{teo:failure-classical-FS} gives
\[
\bigcup_{i=0}^{j}\mathcal S_i
\subset\{M_\alpha f_j>1\}.
\]
Since $|\mathcal S_i|=k^{\lfloor(1-\alpha)i\rfloor}$,
\[
w_*(\mathcal S_i)
=k^{\lfloor(1-\alpha)i\rfloor}k^{-(1-\alpha)i},
\]
and therefore
\begin{equation}\label{eq:wstar-level}
k^{-1}\leq w_*(\mathcal S_i)\leq1.
\end{equation}
It follows that
\[
w_*(\{M_\alpha f_j>1\})
\geq\sum_{i=0}^{j}w_*(\mathcal S_i)
\geq\frac{j+1}{k}.
\]

On the other hand, by \eqref{eq:iterate-wstar} and \eqref{eq:wstar-level},
\begin{align*}
\sum_{x\in T}|f_j(x)|M_\alpha^{(n)}w_*(x)
&\lesssim_{\alpha,k,n}
C_0\sum_{x\in\mathcal S_j}w_*(x)
=C_0w_*(\mathcal S_j)
\lesssim C_0,
\end{align*}
uniformly in $j$. The quotient therefore tends to infinity.
\end{proof}

\begin{remark}\label{rem:extreme-no-s}
The extreme weight $w_*=w_{\alpha-1}$ satisfies $M_\alpha w_*\simeq w_*$, but it does
not satisfy $M_{\alpha,s}w_*\lesssim w_*$ for any $s>1$. For each $j\ge1$, choose
$x_j\in T^j$ and let $x_0$ be its ancestor at distance $j$. The single point
$x_0\in S(x_j,j)$ gives
\begin{align*}
M_{\alpha,s}^{o}w_*(x_j)
\geq |S(x_j,j)|^{\alpha-1/s}w_*(x_0)
\gtrsim_{\alpha,s}
w_*(x_j)k^{j(1-1/s)}.
\end{align*}
Since $1-1/s>0$, the quotient tends to infinity as $j\to\infty$. By
\eqref{eq:Mas-sphere-ball}, the same failure holds for the ball operator
$M_{\alpha,s}$. This is the fractional analogue of the endpoint phenomenon for the
extreme radial weight in \cite[Remark~3.3]{ORRS}.
\end{remark}

\subsection{A necessary condition for weighted weak type $(1,1)$}

\begin{proof}[Proof of Theorem \ref{Teo 1.3}]
Fix $x_0\in T$ and $r\geq0$, and take $f=\delta_{x_0}$. If
$y\in B(x_0,r)$ and $d=d(x_0,y)$, then $d\leq r$ and
\[
M_\alpha(\delta_{x_0})(y)
\geq \frac1{|B(y,d)|^{1-\alpha}}.
\]
Since $|B(y,d)|\leq4k^d\leq4k^r$, we obtain
\[
M_\alpha(\delta_{x_0})(y)
\geq4^{-(1-\alpha)}k^{-(1-\alpha)r}.
\]
Thus, with
\[
\lambda=\frac12\,4^{-(1-\alpha)}k^{-(1-\alpha)r},
\]
we have
\[
B(x_0,r)\subset\{M_\alpha(\delta_{x_0})>\lambda\}.
\]
Applying \eqref{eq:weak-w-assumption},
\[
w(B(x_0,r))
\lesssim_{\alpha,k,C_w}k^{(1-\alpha)r}w(x_0).
\]
Since $|B(x_0,r)|\geq k^r$,
\[
\frac{w(B(x_0,r))}{|B(x_0,r)|^{1-\alpha}}
\lesssim_{\alpha,k,C_w}w(x_0).
\]
Taking the supremum over $r$ proves the result.
\end{proof}

Using  Theorem~\ref{Teo 1.3} and
Theorem~\ref{Teo 1.2} with $n=1$, we obtain the following result.

\begin{corol}\label{cor:necessary-not-sufficient}
If $0<\alpha<1/2$, the condition
\[
M_\alpha w\lesssim w
\]
is necessary but not sufficient for the weighted weak type estimate
\eqref{eq:weak-w-assumption}.
\end{corol}

\subsection{The phase transition at $\alpha=1/2$}

The radial classification suggests that the threshold $\alpha=1/2$ is a
consequence of the geometry of the tree. We now show that no global weight,
radial or otherwise, can satisfy the pointwise condition at or above this
threshold. Above the threshold, the condition forces incompatible rates of
growth. At the endpoint, we use the accumulation of ancestor contributions
in Lemma~\ref{lem:radial-sphere-decomp}.

\begin{prop}\label{prop:no-global-weight}
Let $1/2\le\alpha<1$. There is no nonzero weight $w$ on $T$ satisfying
\[
M_\alpha w(x)\leq Cw(x),\qquad x\in T,
\]
for some finite constant $C$.
\end{prop}

\begin{proof}
First assume that $\alpha>1/2$. Suppose that such a weight exists and choose $x_0$ with $w(x_0)>0$. Let
$D_r(x_0)$ be the set of descendants of $x_0$ at distance $r$; then
$|D_r(x_0)|=k^r$.

For $y\in D_r(x_0)$, since $x_0\in B(y,r)$ and $|B(y,r)|\leq4k^r$,
\[
Cw(y)\geq M_\alpha w(y)
\geq \frac{w(x_0)}{|B(y,r)|^{1-\alpha}}
\gtrsim_{\alpha,k}k^{-(1-\alpha)r}w(x_0).
\]
Summing over $D_r(x_0)$ gives
\begin{equation}\label{eq:lower-growth-weight}
w(B(x_0,r))
\geq w(D_r(x_0))
\gtrsim_{\alpha,k,C}k^{\alpha r}w(x_0).
\end{equation}
On the other hand, applying $M_\alpha w(x_0)\leq Cw(x_0)$ to the ball
$B(x_0,r)$ and using $|B(x_0,r)|\leq4k^r$ yields
\begin{equation}\label{eq:upper-growth-weight}
w(B(x_0,r))
\lesssim_{\alpha,k,C}k^{(1-\alpha)r}w(x_0).
\end{equation}
Combining \eqref{eq:lower-growth-weight} and \eqref{eq:upper-growth-weight} gives
$k^{(2\alpha-1)r}\lesssim1$ for every $r$, a contradiction because
$\alpha>1/2$.

We now consider $\alpha=1/2$. Suppose that
$M_{1/2}w\leq Cw$ for a nonzero weight $w$, and denote the root by $x_0$.
Then $w(x_0)>0$, since some ball centered at $x_0$ contains a vertex where $w$
is positive. For $x\in T^j$, the root belongs to $B(x,j)$, and hence
\[
Cw(x)\geq M_{1/2}w(x)
\geq\frac{w(x_0)}{|B(x,j)|^{1/2}}
\geq\frac{w(x_0)}{2}k^{-j/2}.
\]
Thus, writing $v(x)=k^{-|x|/2}$, we have
\begin{equation}\label{eq:critical-pointwise-lower}
w(x)\geq\frac{w(x_0)}{2C}v(x),\qquad x\in T.
\end{equation}

For $x\in T^j$, formula~\eqref{eq:radial-frac} with
$\alpha=1/2$, $\beta=-1/2$ and $r=j$ gives
\[
\frac1{|S(x,j)|^{1/2}}\sum_{y\in S(x,j)}v(y)
\simeq (j+1)k^{-j/2}.
\]
Since $|B(x,j)|\leq2|S(x,j)|$, it follows that
\[
M_{1/2}v(x)\geq\frac1{|B(x,j)|^{1/2}}
\sum_{y\in S(x,j)}v(y)
\gtrsim (j+1)k^{-j/2}.
\]
Combining this with \eqref{eq:critical-pointwise-lower} and the monotonicity
of the maximal operator, we obtain
\[
Cw(x)\geq M_{1/2}w(x)
\geq\frac{w(x_0)}{2C}M_{1/2}v(x)
\gtrsim\frac{w(x_0)}{C}(j+1)k^{-j/2}.
\]
Summing over the $k^j$ vertices of $T^j$ yields
\begin{equation}\label{eq:critical-generation-lower}
w(T^j)\gtrsim C^{-2}w(x_0)(j+1)k^{j/2}.
\end{equation}
On the other hand, applying the assumed inequality at the root gives
\[
w(T^j)\leq w(B(x_0,j))
\leq Cw(x_0)|B(x_0,j)|^{1/2}
\leq2Cw(x_0)k^{j/2}.
\]
Together with \eqref{eq:critical-generation-lower}, this implies
$j+1\lesssim C^3$ for every $j\geq0$, a contradiction.
\end{proof}

Combining Theorem~\ref{Teo 1.3} with
Proposition~\ref{prop:no-global-weight}, we obtain the following result.

\begin{corol}\label{cor:no-weak-supercritical}
Let $1/2\le\alpha<1$. No nonzero weight can satisfy the weighted weak type
estimate \eqref{eq:weak-w-assumption}.
\end{corol}

\begin{remark}
At the critical value $\alpha=1/2$, the lower and upper exponential rates
in \eqref{eq:lower-growth-weight} and \eqref{eq:upper-growth-weight} coincide.
The additional factor $j+1$ in \eqref{eq:critical-generation-lower} is what
rules out global weights at this endpoint. Together with
Proposition~\ref{prop:radial-weights}, this shows that a nonzero global weight
satisfying $M_\alpha w\lesssim w$ exists if and only if $0<\alpha<1/2$.
By Remark~\ref{rem:nontrivial-weak} and
Corollary~\ref{cor:no-weak-supercritical}, the same range characterizes the
existence of nonzero weights for the weak type estimate
\eqref{eq:weak-w-assumption}.
\end{remark}

\section*{AI Disclosure}
ChatGPT was used for language editing and proofreading and also helped suggest considering subtrees for the example demonstrating the failure of the Fefferman–Stein estimate in Section 3. All mathematical results, arguments, proofs, and conclusions were developed and established by the authors.

\end{document}